\documentclass{article}
\usepackage{amssymb, amsmath}
\newcommand{\pf}{\noindent{\bf Proof.}\ }

\newcommand{\qed}{\begin{flushright}$\Box$\ \ \ \ \ \
	\end{flushright}}

\newcommand{\from}{\leftarrow}
\newcommand{\ybar}{\overline Y}
\newcommand{\unitob}{\mathbf 1}

\newtheorem{thm}{Theorem}[section]
\newtheorem{prop}[thm]{Proposition}
\newtheorem{lemma}[thm]{Lemma}
\newtheorem{cor}[thm]{Corollary}
\newtheorem{dfn}[thm]{Definition}

\usepackage{tikz}
\usetikzlibrary{calc, intersections, patterns}

\usepackage{layout}
\usepackage{dsfont}
\usepackage{enumerate}
\usepackage{graphicx}

\title{{\bf Coisotropic pairs}}
\author{Jonathan Lorand\thanks{Jonathan Lorand was supported by ETH Zurich, the city of Zurich, and the Anna \& Hans K\"agi Foundation. This research was conducted while he was at UC Berkeley as a Visiting Student Researcher.}\\
Department of Mathematics\\
ETH Zurich\\
Zurich, Switzerland\\
\tt{jonathanlorand@gmail.com}
\and 
Alan Weinstein\\Department of Mathematics\\
University of California\\
Berkeley, CA 94720 USA\\ \tt{alanw@math.berkeley.edu}
}
\date{}
\begin{document}
\maketitle
\begin{abstract}
We give two equivalent sets of invariants which classify pairs of coisotropic subspaces of a finite-dimensional symplectic vector space. We identify five elementary types of coisotropic pairs and show that any coisotropic pair decomposes in an appropriate sense as the direct sum of coisotropic pairs of elementary type. 
\end{abstract}

\section{Introduction}

The classification given in this short note is a beginning step in two separate projects, not yet complete.

The first project is a classification, up to conjugation by linear symplectomorphisms, of canonical relations (lagrangian correspondences) from a finite-dimensional symplectic vector space to itself.   Without symplectic structure, this
classification of linear relations was carried out by Towber \cite{Towber}.  In the symplectic situation, for the special case of graphs of symplectomorphisms, the classification amounts to identifying the conjugacy classes in the group of symplectic matrices. This classification and the problem of finding associated normal forms has a long history extending from Williamson \cite{Williamson} to Gutt \cite{Gutt}. In the general symplectic case, a result of Benenti and Tulczyjew (\cite{Benenti Tulczyjew}, Proposizioni 4.4 \& 4.5) tells us that a canonical relation $X\from Y$ is given by coisotropic subspaces of $X$ and $Y$ and a symplectomorphism between the corresponding reduced spaces.  When $X = Y,$ a first step in the classification of canonical relations is then a classification of the coisotropic pairs giving the range and domain.

The second project is an extension of the Wehrheim-Woodward theory of linear canonical relations (see \cite{Li-Bland Weinstein}, \cite{Weinstein}) to the case where the set of lagrangian correspondences $X\from Y$ is replaced by the set of coisotropic correspondences, i.e. coisotropic subspaces of
$X \times \ybar.$ Each pair of coisotropic subspaces of $X$ gives a WW morphism represented by a diagram of the form
$\unitob \from X \from \unitob$, and isomorphic pairs correspond to the same WW morphism. There are also inequivalent pairs representing the same WW morphism. The problem is to determine exactly which pairs are ``WW equivalent''. 

In the following we consider the situation when $X = Y$ and denote this finite dimensional vector space by $V,$ equipped with symplectic form $\omega$. An ordered pair $(A,B)$ of coisotropic subspaces $A$ and $B$ in $V$ will be called a \textbf{\emph{coisotropic pair}}. Coisotropic pairs $(A,B)$ and $(\hat A, \hat B)$ given in $(V,\omega)$ and $(\hat V, \hat \omega)$ respectively are \textbf{\emph{equivalent}} if there exists a linear symplectomorphism $S:V \rightarrow \hat V$ such that $S(A) = \hat A$ and $S(B) = \hat B$. We will show that a coisotropic pair $(A,B)$ in $(V,\omega)$ is fully characterized up to this equivalence by the following five numbers
$$\dim (A^{\omega} \cap B^{\omega}), \ \ \dim A^{\omega}, \ \  \dim B^{\omega}, \ \  \frac{1}{2}\dim V, \ \ \dim (A^{\omega} \cap B)$$ 
where for any linear subspace $W \subset V$, $ W^{\omega}$ denotes its (symplectic) \textbf{\emph{orthogonal}} $\{v \in V \mid \omega(v,w) = 0 \ \forall w \in W\} $.  
We call these five numbers the \textbf{\emph{canonical invariants}} of $(A,B)$ and label them $k_1$ through $k_5$ in the above order. They are largely independent, subject only to certain inequalities (see Corollary \ref{ineq only constraint}). 
 
The first four invariants $k_1, k_2, k_3, k_4$ characterize the subspaces $A$ and $B$ up to the above equivalence if one drops the condition that $S$ be symplectic and that $A$ and $B$ be coisotropic, i.e. these four invariants contain the purely linear algebraic information.  Indeed, using the identities $\dim W^{\omega} = \dim V - \dim W$ and $(E + F)^{\omega} = E^{\omega} \cap F^{\omega}$, which hold for any subspaces $W$, $E$, $F \subset V$, one can obtain the the linear algebraic data 
\begin{equation}\label{lin alg invariants}
\dim V, \ \dim A, \ \dim B, \ \dim (A \cap B)
\end{equation}
from these four invariants:
$\dim V = 2 \cdot \frac{1}{2} \dim V$, $\dim A = \dim V - \dim A^{\omega}$, $\dim B = \dim V - \dim B^{\omega}$,
and
\begin{align*}
\dim (A \cap B) & = \dim V - \dim (A \cap B)^{\omega} \\
& = \dim V - \dim (A^{\omega} + B^{\omega}) \\
 & = \dim V - \dim A^{\omega} - \dim B^{\omega} + \dim (A^{\omega} \cap B^{\omega})
\end{align*}
It is straightforward to check that his relationship is invertible; one could thus equivalently use the numbers (\ref{lin alg invariants}) as the first four invariants. 

The fifth invariant $k_5 = \dim (A^{\omega} \cap B)$ is what fixes the symplectic information. One could equivalently choose $ \dim (B^{\omega} \cap A)$ as the fifth invariant, since
\begin{align*}
 \dim (A^{\omega} \cap B) & = \dim V - \dim (A^{\omega} \cap B)^{\omega} \\
  & = \dim V - [\dim A + \dim B^{\omega} -  \dim (B^{\omega} \cap A)]  \\
  & = \dim (B^{\omega} \cap A) + \dim A^{\omega} -\dim B^{\omega}
\end{align*} 

When $\dim A = \dim B$ \footnote{This is the case, for example, when $A$ and $B$ are the range and domain respectively of a linear canonical relation in $V \times \bar V.$}, it follows that $\dim (A^{\omega} \cap B) = \dim (B^{\omega} \cap A)$, and a total of four invariants suffice to characterize the coisotropics $A$ and $B$. They can be given in a symmetric way as 
$$
\dim V, \ \dim (A + B), \ \dim (A \cap B), \ \text{rank} (A^{\omega} + B^{\omega})
$$
where for any subspace $W \subset V$, rank$(W) = \dim W - \dim (W \cap W^{\omega})$. The symmetry of these invariants implies that $(A,B)$ and $(B,A)$ are equivalent as coisotropic pairs when $\dim A = \dim B$. 

Because a coisotropic subspace $A$ is uniquely determined by the isotropic subspace $A^{\omega}$, and $S(A^{\omega}) = S(A)^{\hat \omega}$ for any linear symplectomorphism $S:V \rightarrow \hat V$, one could equivalently consider isotropic pairs instead of coisotropic ones. This indeed simplifies some calculations and proofs; for the present though we treat things from the coisotropic standpoint.

We think in terms of the ground field $\mathbb{R}$, though all results should hold for other fields, with the exception perhaps of characteristic 2. One may also include, under slight modifications, the situation where the symplectic form is replaced by a non-degenerate symmetric bilinear pairing. 

For convenience, all maps and subspaces are tacitly assumed linear unless otherwise stated, and a linear symplectomorphism will be synonymously called a symplectic map. The letters $A$ and $B$ always denote coisotropic subspaces of a finite-dimensional symplectic vector space $(V,\omega)$. Subspaces $E,F \subset V$ are called \textbf{$\omega$-\emph{orthogonal}} if $\omega(e,f) = 0 \ \forall e \in E \ \forall f \in F$. We use the notation $\tilde \omega$ for the isomorphism $V \rightarrow V^*, \ v \mapsto \omega (v, \cdot) $ induced by the symplectic form $\omega$, and the symbol `$\simeq$' denotes a linear isomorphism, not necessarily symplectic.

\section{Classification of coisotropic pairs}

We recall quickly some basic facts from symplectic linear algebra which will be useful in the following. Let $W, E$ and $F$ be subspaces of $V$ and denote by $\omega_W$ the restriction of $\omega$ to $W$. One has $\ker \omega_W = W \cap W^{\omega} $, and the reduced space $W/(W \cap W^{\omega})$ admits an induced symplectic form $[\omega]$ given by 
$$[\omega]([u], [v] ) := \omega (u,v) \quad  \quad u,v \in W$$

In addition to symplectic reduction $W \rightarrow W/(W\cap W^{\omega})$, a different quotient relationship arises from the isomorphism $\tilde \omega : V \rightarrow V^*.$ Namely, $\tilde \omega$ post-composed with the restriction $V^* \rightarrow W^*$ has kernel $W^{\omega}$, hence it induces a natural isomorphism $V/W^{\omega} \rightarrow W^*.$ In the special case when $W$ is a lagrangian subspace, $V/W \simeq W^*$. If $(L,L')$ is a transversal lagrangian pair in V, i.e. lagrangian subspaces such that $V= L \oplus L'$, then $V \simeq L \oplus L^*$ symplectically via the natural map 
$$ id \oplus \tilde \omega : V  \rightarrow L \oplus L^*, \quad u + v \mapsto (u,  \omega (v, \cdot)) $$
where the space $L \oplus L^*$ is endowed with the symplectic form 
$$ ((v, \alpha), (w, \beta)) \mapsto \beta (v) - \alpha (w) $$ 
With respect to this form, $L\times 0$ and $0 \times L^*$ are lagrangian subspaces and they are the images of $L$ and $L'$ respectively under the map $id \oplus \tilde \omega$ above. In particular it follows that for any two transversal lagrangian pairs $(L,L')$ and $(\hat L, \hat L')$ given  in symplectic spaces $V$ and $\hat V$ of the same dimension there always exists a symplectic map $S: V \rightarrow \hat V$ such that $S(L) = \hat L$ and $S(L') = \hat L'$. 

In general, if $V = E \oplus F$ and $\hat V = \hat E \oplus \hat F$, we say that a map $S:V \rightarrow \hat V$ satisfying $S(E) = \hat E$ and $S(F) = \hat F$ \textbf{\emph{respects the decompositions}} in $V$ and $\hat V$. If $E$, $F$, $\hat E$ and $\hat F$ are symplectic and $S$ is a symplectic map which respects the decompositions, then $S\vert_E: E \rightarrow \hat E$ and $S \vert_F : F \rightarrow \hat F$ are symplectic maps. On the other hand, if $E,F$ are $\omega$-orthogonal, $\hat E,\hat F$ $\hat \omega$-orthogonal, and $\sigma : E \rightarrow \hat E$, $\rho: F \rightarrow \hat F$ are symplectic maps, then $\sigma \oplus \rho$ defines a symplectic map $V \rightarrow \hat V$ which respects the decompositions in $V$ and $\hat V$. The $\omega$-orthogonality condition on $E$ and $F$ (and $\hat E$ and $\hat F$) amounts to $E \oplus F$ being naturally symplectomorphic to the external direct sum of two separate symplectic spaces $(E,\omega \vert_E)$ and $(F, \omega \vert_F)$, endowed with the direct sum symplectic form $\omega \vert_E \oplus \omega \vert_F$ defined by $((e,f),(e',f')) \mapsto \omega_E(e,e') + \omega_F(f,f')$.

A useful way to obtain $\omega$-orthogonal direct sum decompositions is the following.\footnote{see \cite{Cushman Bates}, p. 275 and p. 404, as well as \cite{Artin}, p. 120 and \cite{Witt}, p. 33 (Satz 1).} 

\begin{lemma}[Witt-Artin decomposition]\label{Witt-Artin} 
Let $W \subset V$ be any subspace, and $E$ and $F$ complements of $W\cap W^{\omega}$ in $W$ and $ W^{\omega}$ respectively. Then $E$ and $F$ are symplectic subspaces and $\omega$-orthogonal, and $V$ decomposes as the $\omega$-orthogonal direct sum 
$$ V = E \overset{\omega}{\oplus} F \overset{\omega}{\oplus} (E \oplus F)^{\omega}$$ 
Moreover, $W\cap W^{\omega}$ is a lagrangian subspace of $(E \oplus F)^{\omega}.$ 
\end{lemma}

\pf
Let $\pi: W \rightarrow E$ be the projection map associated to the decomposition $W = W\cap W^{\omega} \oplus E$. This induces an isomorphism $\tilde \pi : W/W\cap W^{\omega} \rightarrow E$ such that $\tilde \pi ([v]) = \pi(v)$ for all $v\in W$. Under this map, the symplectic form $[\omega]$ on the reduced space $W/W\cap W^{\omega}$ is pushed forward to a symplectic form on $E$, and $\tilde \pi_*[\omega] = \omega_E$:  
$$\tilde \pi_*[\omega](e_1,e_2) = [\omega](\tilde \pi^{-1}e_1,\tilde \pi^{-1}e_2) = [\omega]([e_1],[e_2]) = \omega(e_1,e_2) \quad \quad \forall e_1,e_2 \in E $$
Thus $E$ is symplectic, and by analogous arguments $F$ is symplectic as well. Because $E \subset W$ and $F \subset W^{\omega}$,  $E$ and $F$ are $\omega$-orthogonal. As a consequence, $E \cap F = 0$ and $E \oplus F$ is symplectic also. From this it follows that $V = E \oplus F \oplus (E \oplus F)^{\omega}$. 

Finally, $W\cap W^{\omega}$ is in $(E \oplus F)^{\omega}$ since it is in $E^{\omega}$ and $F^{\omega}$ each, and $(E + F)^{\omega} = E^{\omega} \cap F^{\omega}$. Clearly $W\cap W^{\omega}$ is isotropic. To see that it is lagrangian in $G:=(E \oplus F)^{\omega}$, note that $W^{\omega} + W = E \oplus F \oplus (W \cap W^{\omega})$ and recall the general fact that if $U,X,Y$ are subspaces such that $U \supset X$ and $X \cap Y = 0$,  it holds that $U \cap (X \oplus Y) = X \oplus (U \cap Y)$. We now calculate 
\begin{align*}
(W\cap W^{\omega})^{\omega_{_G}} & = (W\cap W^{\omega})^{\omega} \cap G \\ 
& = [(E \oplus F) \oplus (W \cap W^{\omega})] \cap G \\ 
& = [(E \oplus F) \cap G] \oplus [W \cap W^{\omega}] \\ 
& = (W\cap W^{\omega})
\end{align*}
where the last inequality uses the fact that $(E\oplus F) \cap (E \oplus F)^{\omega} = 0$ and the second to last uses the general fact about subspaces above, with $G$ in the role of $U$. 
\qed

\begin{center}
\begin{tikzpicture}[scale = 1]

\path[name path=V, draw] (2.5,2.5) circle [radius=2];

\path[name path=V, clip] (2.5,2.5) circle [radius=2];

\coordinate (c) at (2.5,2.5);
\coordinate (r1) at (0,2.5);
\coordinate (r2) at (5,2.5);
\coordinate (r3) at (0,4.5);
\coordinate (r4) at (0,.5);
\coordinate (r5) at (5,4);
\coordinate (r6) at (5,1);
\coordinate (r7) at (1.5,5);
\coordinate (r8) at (1.5,0);
\coordinate (r9) at (3,5);
\coordinate (r10) at (3,0);

\fill[color=red!90!white!, opacity=.6] (c) -- (r3) -- (r4) -- (r10) -- (r6) -- cycle;
\fill[color=yellow!90!white!, opacity=.6] (c) -- (r5) -- (r9) -- (r3) -- cycle;
\fill[color=blue!90!white!, opacity=.6] (c) -- (r5) -- (r6) -- (r6) -- (r10) -- (r4) -- cycle;

\path[name path=V, draw] (2.5,2.5) circle [radius=2];
\draw (c) -- (r3);
\draw (c) -- (r4);
\draw (c) -- (r5);
\draw (c) -- (r6);

\node at (2.6,1.2) {\small{$W \cap W^{\omega}$}};
\node at (2.5,3.8) {\small{$L'$}};
\node at (1.1,2.5) {\small{$E$}};
\node at (3.9,2.5) {\small{$F$}};

\end{tikzpicture}

Figure 1

\end{center}

The Witt-Artin decomposition of $V$ with respect to $W$ is represented diagrammatically in Figure 1, each piece representing a direct summand. The circle is all of $V,$ red is used for the subspace $W$ and blue for $W^{\omega},$ giving a violet hue where they intersect. The yellow subspace $L'$ represents a choice of a lagrangian complement of $W\cap W^{\omega}$ in $(E \oplus F)^{\omega}$.

We now turn to our objects of study, two coisotropic subspaces $A$ and $B$ of $V$, fixing the notation $I := A^{\omega} \cap B^{\omega} $ and $K:= A^{\omega} \cap B + B^{\omega} \cap A$. As announced, the numbers $\dim (A^{\omega} \cap B^{\omega})$, $\dim A^{\omega}$, $\dim B^{\omega}$, $ 1/2 \dim V$ and $\dim (A^{\omega} \cap B)$, which we call the canonical invariants associated to $(A,B)$, completely characterize a coisotropic pair up to equivalence.

\begin{prop}\label{main prop}
Let $(A,B)$ and $(\hat A,\hat B)$ be pairs of coisotropic subspaces in $(V,\omega)$ and $(\hat V, \hat \omega)$ respectively. Then $(A,B)$ and $(\hat A,\hat B)$ are equivalent if and only if their associated canonical invariants are equal. 
\end{prop}

\pf
If $(A,B)$ and $(\hat A,\hat B)$ are equivalent via some symplectic map $S: V \rightarrow \hat V$, it is clear that all the canonical invariants of $(A,B)$ and $(\hat A,\hat B)$ coincide. 

For the converse, we will show that $V$ can be written as an $\omega$-orthogonal direct sum of five symplectic subspaces
$$ V = D \overset{\omega}{\oplus} E \overset{\omega}{\oplus} F \overset{\omega}{\oplus} G \overset{\omega}{\oplus} H $$
where each symplectic piece, except for $F$, is further decomposed as a lagrangian pair
$$  D = I \oplus J, \ \ \  E = E_1 \oplus E_2 , \ \ \  G = G_1 \oplus G_2, \ \ \  H = H_1 \oplus H_2 $$
so that we obtain a decomposition of $V$ into a total of nine subspaces
\begin{equation}\label{inner decomp}
V = (I \oplus J) \overset{\omega}{\oplus} (E_1 \oplus E_2) \overset{\omega}{\oplus} F \overset{\omega}{\oplus} (G_1 \oplus G_2) \overset{\omega}{\oplus} (H_1 \oplus H_2)
\end{equation}
Moreover, this decomposition will have the following properties:

\begin{enumerate}[i)]
\item  the dimension of each summand is uniquely determined by the canonical invariants of $(A,B)$
\item $A$ and $B$ are decomposable as
$$A = I \oplus E_1 \oplus G_1 \oplus F \oplus H_1 \oplus H_2$$ 
$$ B = I \oplus E_2  \oplus H_1 \oplus F  \oplus G_1 \oplus G_2$$ 

\end{enumerate} 

One can decompose $\hat V$ in an analogous manner, and hence when $(A,B)$ and $(\hat A ,\hat B)$ have the same invariants, by property i) the dimensions of corresponding symplectic pieces in the decompositions of $V$ and $\hat V$ will match. In this case, for dimension reasons alone there exist five symplectic maps, one each between corresponding symplectic pieces, i.e. one from $D$ to $\hat D$, one from $E$ to $\hat E$, and so on.  These maps can further be be chosen to respect the respective decompositions into lagrangian pairs. 

Because the five-part decompositions of $V$ and $\hat V$ are $\omega$-orthogonal, the direct sum of these five symplectic maps defines a symplectic map $S: V \rightarrow \hat V$ which respects all nine summands of the decompositions of $V$ and $\hat V$. In particular, by property ii), $S$ will then also satisfy $S(A)=\hat A$ and $S(B)=\hat B$.

To achieve the decomposition (\ref{inner decomp}) we will construct a certain Witt-Artin decomposition of $V$ with respect to $W := A^{\omega} + B^{\omega}$, refined and adapted to the coisotropic subspaces $A$ and $B$.

Recall that $I = A^{\omega} \cap B^{\omega} $ and $K= A^{\omega} \cap B + B^{\omega} \cap A$, and note that $$W^{\omega} = (A^{\omega} + B^{\omega})^{\omega} = A \cap B$$ 
and
$$W \cap  W^{\omega} = (A^{\omega} + B^{\omega}) \cap (A \cap B) = A^{\omega} \cap B + B^{\omega} \cap A = K$$ 
We begin by decomposing $A^{\omega}$ into three parts by choosing a subspace $G_1$ such that $A^{\omega} \cap B = I \oplus G_1$ and a subspace $E_1$ such that $A^{\omega} = A^{\omega} \cap B \oplus E_1$, giving a decomposition 
$$ A^{\omega} = I \oplus G_1 \oplus E_1 $$
Analogously we obtain a decomposition 
$$B^{\omega} = I \oplus H_1 \oplus E_2$$
where $E_2$ is such that $B^{\omega} = B^{\omega} \cap A \oplus E_2$, and $H_1$ such that $B^{\omega} \cap A = I \oplus H_1$. Note that $H_1$ and $G_1$ have zero intersection, since $H_1 \cap G_1 \subset A^{\omega} \cap B^{\omega} = I$ and $H_1 \cap I = 0$ and $G_1\cap I = 0$. Similarly, $E_1 \cap E_2 = 0$. In particular we have 
$$K = A^{\omega} \cap B + B^{\omega} \cap A =  I \oplus G_1 + I \oplus H_1 =  I \oplus G_1 \oplus H_1$$ 

We now set $E:= E_1 \oplus E_2$. This defines a subspace such that $K \oplus E = A^{\omega} + B^{\omega} = W$. Indeed, 
$$K + E = I \oplus G_1 \oplus H_1 + E_1 \oplus E_2 = A^{\omega} + B^{\omega} = W$$ and $K \cap E = 0$ since 
\begin{align*}
\dim K + \dim E & = (\dim I  + \dim G_1 + \dim H_1) + (\dim E_1 + \dim E_2) \\
& = \dim A^{\omega} + \dim B^{\omega} - \dim I \\
& = \dim (A^{\omega} + B^{\omega}) = \dim W \\
& = \dim (K + E) 
\end{align*}
Because $E$ is a complement of $K =W\cap W^{\omega}$ in $W$, $E$ is symplectic by Lemma \ref{Witt-Artin}, and since $E_1$ and $E_2$ are both isotropic, we conclude that they form a transversal lagrangian pair in $E$. 

To obtain a Witt-Artin decomposition with respect to $W$, we choose a complement $F$ of $W \cap W^{\omega} = K$ in $W^{\omega} = A \cap B$, i.e. so that 
$$ A \cap B  =  K \oplus F $$
Applying Lemma \ref{Witt-Artin} again we know that $F$ is symplectic, as is $E$, and $V$ decomposes into the $\omega$-orthogonal direct sum
$$ V = E \overset{\omega}{\oplus} F  \overset{\omega}{\oplus} (E \oplus F)^{\omega} $$
with $K$ as a lagrangian subspace of the symplectic subspace $(E \oplus F)^{\omega}$.

We refine this decomposition by choosing a lagrangian complement $K'$ of $K$ in $(E \oplus F)^{\omega}$ and by defining a decomposition in $K'$ using the decomposition $K =  I \oplus G_1 \oplus H_1$ as follows. Any basis $\textbf{q}$ of $K$ is mapped under $\tilde{\omega}$ to a basis of $(K')^*$, whose dual basis $\textbf{p}$ in $K'$ is conjugate to $\textbf{q}$, \text{i.e.} together $\textbf{q}$ and $\textbf{p}$ form a symplectic basis of $K \oplus K'$.  If we consider a basis $\textbf{q}$ which is adapted to the decomposition in $K$, then this partitioning induces a partitioning of $\textbf{p}$ which defines subspaces $J$, $G_2$ and $H_2$ in $K'$ such that 
$$
K' = J \oplus G_2 \oplus H_2
$$
and $D:= I \oplus J$, $G:= G_1 \oplus G_2$ and $H:=H_1\oplus H_2$ are $\omega$-orthogonal symplectic subspaces, comprised each of a lagrangian pair, giving
$$
K \oplus K' = D \oplus G \oplus H
$$
In total we thus obtain a decomposition 
$$ V = (I \oplus J ) \oplus (E_1 \oplus E_1) \oplus F \oplus (G_1 \oplus G_2) \oplus (H_1 \oplus H_2) $$ 
where parentheses enclose transversal lagrangian pairs in a symplectic subspace. This decomposition is visualized in Figure 2 - the full circle represents $V$, each piece is a direct summand, and lagrangian pairs are aligned symmetrically with respect to the horizontal axis and shaded with colors of a similar hue. 

\begin{center}
\begin{tikzpicture}[scale = 1]

\path[name path=V, draw] (2.5,2.5) circle [radius=2];

\path[name path=V, clip] (2.5,2.5) circle [radius=2];

\coordinate (c) at (2.5,2.5);
\coordinate (r1) at (0,2.5);
\coordinate (r2) at (5,2.5);
\coordinate (r3) at (0,4.5);
\coordinate (r4) at (0,.5);
\coordinate (r5) at (5,4);
\coordinate (r6) at (5,1);
\coordinate (r7) at (1.5,5);
\coordinate (r8) at (1.5,0);
\coordinate (r9) at (3,5);
\coordinate (r10) at (3,0);

\fill[color=brown!80!white!, opacity=.4] (c) -- (r3) -- (r1) -- cycle; 
\fill[color=brown!80!gray!, opacity=.6] (c) -- (r4) -- (r1) -- cycle; 
\fill[color=yellow!80!white!, opacity=.4] (c) -- (r7) -- (r3) -- cycle; 
\fill[color=yellow!80!gray!, opacity=.6] (c) -- (r8) -- (r4) -- cycle; 
\fill[color=orange!80!white!, opacity=.4] (c) -- (r9) -- (r7) -- cycle; 
\fill[color=orange!80!gray!, opacity=.6] (c) -- (r10) -- (r8) -- cycle; 
\fill[color=red!80!white!, opacity=.4] (c) -- (r5) -- (r9) -- cycle; 
\fill[color=red!80!gray!, opacity=.6] (c) -- (r6) -- (r10) -- cycle; 
\fill[color=violet!80!brown!, opacity=.5] (c) -- (r5) -- (r6) -- cycle; 

\path[name path=V, draw] (2.5,2.5) circle [radius=2];
\draw (c) -- (r1);
\draw (c) -- (r3);
\draw (c) -- (r4);
\draw (c) -- (r5);
\draw (c) -- (r6);
\draw (c) -- (r7);
\draw (c) -- (r8);
\draw (c) -- (r9);
\draw (c) -- (r10);

\node at (2.35,1.2) {\small{$I$}};
\node at (2.35,3.8) {\small{$J$}};
\node at (1.2,2) {\small{$E_1$}};
\node at (1.2,2.9) {\small{$E_2$}};
\node at (1.7,1.4) {\small{$G_1$}};
\node at (3.3,1.3) {\small{$H_1$}};
\node at (3.9,2.5) {\small{$F$}};
\node at (1.7,3.6) {\small{$G_2$}};
\node at (3.3,3.6) {\small{$H_2$}};

\end{tikzpicture}

Figure 2

\end{center}
\ 

The coisotropics $A$ and $B$ are related to the decomposition in $K'$ in that $G_2 = B \cap K'$ and $H_2 = A \cap K'$. To see this it suffices to show the corresponding equalities for the orthogonal spaces. For the case of $A\cap K'$ (the case for $B \cap K'$ is analogous) one has
\begin{align*}
(A \cap K' )^{\omega} & = A^{\omega} + (K')^{\omega} \\
				& = I \oplus G_1 \oplus E_1 + E \oplus F \oplus K' \\
				& = I \oplus G_1 \oplus E \oplus F \oplus K' \\
				& = H_2^{\omega}
\end{align*}
where we use in the last step that $H_2$ is $\omega$-orthogonal to $D$, $G$, $K'$ and $E \oplus F$ and that the dimensions match. 
\

It can now be quickly checked that our decomposition of $V$ satisfies property ii), \text{i.e.} that
$$A = I \oplus E_1 \oplus G_1 \oplus F \oplus H_1 \oplus H_2$$ 
$$ B = I \oplus E_2  \oplus H_1 \oplus F  \oplus G_1 \oplus G_2$$ 
We show this for $A$, the decomposition of $B$ follows in the same way. The inclusion ``$\supset$" is obvious since all the spaces on the right-hand side are subsets of $A$. The opposite inclusion ``$\subset$" can be argued using dimensions:
\begin{align*}
\dim A & =  \dim V - \dim A^{\omega} \\
	       & =  \dim(I \oplus E_1 \oplus G_1  \oplus F \oplus H_1 \oplus H_2) + \dim (J \oplus G_2 \oplus E_2) - \dim A^{\omega} \\
	        & =  \dim(I  \oplus E_1 \oplus G_1  \oplus F \oplus H_1 \oplus H_2) 
\end{align*}
where the last equality follows from the fact that 
$$ \dim A^{\omega} = \dim (I  \oplus E_1 \oplus G_1 ) = \dim (J \oplus E_2 \oplus G_2) $$
since $\dim I = \dim J$, $ \dim E_1 = \dim E_2$, and $\dim G_1 = \dim G_2$ (each pair of subspaces is a lagrangian pair in $D$, $E$ and $G$ respectively). 

The decompositions of $A$ and $B$ are visualized below.  
\begin{center}
\begin{minipage}{.45\textwidth}
\begin{center}

\begin{tikzpicture}[scale = .95]

\path[name path=V, draw] (2.5,2.5) circle [radius=2];

\path[name path=V, clip] (2.5,2.5) circle [radius=2];

\coordinate (c) at (2.5,2.5);
\coordinate (r1) at (0,2.5);
\coordinate (r2) at (5,2.5);
\coordinate (r3) at (0,4.5);
\coordinate (r4) at (0,.5);
\coordinate (r5) at (5,4);
\coordinate (r6) at (5,1);
\coordinate (r7) at (1.5,5);
\coordinate (r8) at (1.5,0);
\coordinate (r9) at (3,5);
\coordinate (r10) at (3,0);

\fill[color = yellow!70!white!, opacity = .3] (2.5,2.5) circle [radius=2]; 

\fill[color=blue!50!white!, opacity=.5] (c) -- (r4) -- (0,0) -- (5,0) -- (r5) -- cycle; 
\fill[color=blue!50!white!, opacity=.5] (c) -- (r7) -- (r3) -- cycle; 
\fill[color=blue!50!white!, opacity=.5] (c) -- (r3) -- (r1) -- cycle; 
\fill[color=green!50!white!, opacity=.5] (c) -- (r1) -- (0,0) -- (5,0) -- (5,5) -- (3,5) -- cycle; 

%
%
%

\fill[color=cyan!70!white!, opacity=.3] (c) -- (r3) -- (r1) -- cycle; 
\fill[color=green!70!yellow!, opacity=.2] (c) -- (r4) -- (r1) -- cycle; 

\fill[color=green!70!yellow!, opacity=.2] (c) -- (r8) -- (r4) -- cycle; 
\fill[color=green!70!yellow!, opacity=.2] (c) -- (r10) -- (r8) -- cycle; 
\fill[color=cyan!70!white!, opacity=.3] (c) -- (r10) -- (r8) -- cycle; 
\fill[color=cyan!70!white!, opacity=.3] (c) -- (r6) -- (r10) -- cycle; 




\path[name path=V, draw] (2.5,2.5) circle [radius=2];
\draw (c) -- (r1);
\draw (c) -- (r3);
\draw (c) -- (r4);
\draw (c) -- (r5);
\draw (c) -- (r6);
\draw (c) -- (r7);
\draw (c) -- (r8);
\draw (c) -- (r9);
\draw (c) -- (r10);

\node at (2.35,1.2) {\small{$I$}};
\node at (2.35,3.8) {\small{$J$}};
\node at (1.2,2) {\small{$E_1$}};
\node at (1.2,2.9) {\small{$E_2$}};
\node at (1.7,1.4) {\small{$G_1$}};
\node at (3.3,1.3) {\small{$H_1$}};
\node at (3.9,2.5) {\small{$F$}};
\node at (1.7,3.6) {\small{$G_2$}};
\node at (3.3,3.6) {\small{$H_2$}};

\end{tikzpicture}

Figure 3

\end{center}
\end{minipage}%
\begin{minipage}{.45\textwidth}
\begin{center}

\begin{tikzpicture}[scale = .7]

\path[use as bounding box] (0,0) rectangle (8, 5.2);

\draw[name path=border] (0,0) rectangle (8,5.2);
\fill [name path=border, color=yellow!70!white!, opacity=.3] (0,0) rectangle (8,5.2);

\path[name path global=B] (3.5,2.5) ellipse [x radius=3, y radius=1, rotate=40]; 
\path[name path global=A] (4,2.5) ellipse [x radius=4, y radius=1.5, rotate=155];
\path[name path global=B'] (4,2) ellipse [x radius=4, y radius=1, rotate=40]; 
\path[name path global=A'] (3.5,2) ellipse [x radius=5, y radius=1.5, rotate=155]; 


\path[name path global=B, fill] [color=blue!50!white!, opacity=.5] (3.5,2.5) ellipse [x radius=3, y radius=1, rotate=40]; 
\path[name path global=A, fill] [color=green!50!white!, opacity=.5] (4,2.5) ellipse [x radius=4, y radius=1.5, rotate=155];
\begin{scope}
\path[name path global=B, clip] (3.5,2.5) ellipse [x radius=3, y radius=1, rotate=40];
\fill [name intersections={of=A' and B, name=q, total=\t}, color=cyan!70!white!, opacity=.3]
(q-4) -- (3,5.2) -- (8,5.2) -- (8,3) -- (q-3) -- cycle;
\end{scope}

\begin{scope}
\path[name path global=A, clip] (4,2.5) ellipse [x radius=4, y radius=1.5, rotate=155]; 
\fill [name intersections={of=A and B', name=p, total=\t}, color=green!70!yellow!, opacity=.2]
(p-1) -- (0,2) -- (0,5.2) -- (5,5.2) -- (p-3) -- cycle;
\end{scope}

\path[name path global=B, draw] (3.5,2.5) ellipse [x radius=3, y radius=1, rotate=40]; 
\path[name path global=A, draw] (4,2.5) ellipse [x radius=4, y radius=1.5, rotate=155];

\draw [name intersections={of=A' and B, name=q, total=\t}]
(q-4) -- (q-3);
\draw [name intersections={of=A and B', name=p, total=\t}]
(p-1) -- (p-3);

\node at (4.1,3.7) {\small{$I$}};
\node at (7,3.5) {\small{$J$}};
\node at (1.7,3.8) {\small{$E_1$}};
\node at (5.1,4.1) {\small{$E_2$}};
\node at (2.9,2.8) {\small{$G_1$}};
\node at (4.9,3.3) {\small{$H_1$}};
\node at (3.8,2.3) {\small{$F$}};
\node at (1.9,1.1) {\small{$G_2$}};
\node at (6,1.5) {\small{$H_2$}};


%
%
%
%

\end{tikzpicture}

Figure 4

\end{center}
\end{minipage}
\end{center}

Figure 3 is a recoloring of Figure 2, and Figure 4 gives an intuitive representation of $A$ and $B$ intersecting, where $V$ is given by the entire rectangle. This is not a proper Venn diagram in the set-theoretic sense, though certain intersections are represented properly, namely $A^{\omega} \cap B$, $B^{\omega} \cap A$ and $A^{\omega} \cap B^{\omega} = I$. 

It remains now only to check that the property i) is fulfilled, \text{i.e.} that the dimensions of the nine summands in our decomposition are uniquely determined by the canonical invariants associated to the pair $(A,B)$. Since any lagrangian subspace of a symplectic subspace has half the dimension of the space within which it is lagrangian, it suffices to show for example that the dimensions of the  subspaces $I$, $E$, $F$, $G_1$ and $H_1$ are uniquely determined.

First, 
$$\dim I = \dim (A^{\omega} \cap B^{\omega}) = k_1$$ 
and the relationships  
\begin{align*}
\dim K & = \dim (A^{\omega} \cap B) + \dim (B^{\omega} \cap A) - \dim I \\
& =  \dim (A^{\omega} \cap B) + [\dim (A^{\omega} \cap B)+ \dim B^{\omega} - \dim A^{\omega}] - \dim I \\
& = 2k_5 + k_3 - k_2 - k_1
\end{align*}
and 
\begin{align*}
\dim W & =  \dim (A^{\omega} + B^{\omega}) \\
& = \dim A^{\omega} + \dim B^{\omega} - \dim  (A^{\omega} \cap B^{\omega}) \\
& = k_2 + k_3 - k_1
\end{align*}
show that $\dim K$ and $\dim W$ are determined.  

Because $E \simeq W/K$ and $F \simeq W^{\omega}/K$ we have 
$$
\dim E = \dim W - \dim K = 2k_2 - 2k_5 
$$
and 
$$
\dim F = \dim (A \cap B) - \dim K = 2k_1 - 2k_3 + 2k_4 - 2k_5
$$

Lastly, $G_1 \simeq (A^{\omega} \cap B) / I$ and $H_1 \simeq (B^{\omega} \cap A) / I$, so $$\dim G_1 = \dim (A^{\omega} \cap B) - \dim I = k_5 - k_1$$ 
and 
$$\dim H_1 = \dim (B^{\omega} \cap A) - \dim I = - k_1 - k_2 + k_3 + k_5$$ 
which proves the property i) and concludes the proof. \qed

\section{Elementary types and normal forms}

The key to Proposition \ref{main prop} was the decomposition (\ref{inner decomp}), satisfying the properties i) and ii). One may rephrase the construction as follows. We found an $\omega$-orthogonal decomposition 
$$V = V_1 \oplus V_2 \oplus V_3 \oplus V_4 \oplus V_5$$ 
into five symplectic subspaces, such that 

a) the dimensions of these subspaces are uniquely determined by the canonical invariants associated to the coisotropic pair $(A,B)$, and 

b) $A$ and $B$ decompose into direct sums
$$ A = A_1 \oplus A_2 \oplus A_3 \oplus A_4 \oplus A_5 $$
$$ B = B_1 \oplus B_2 \oplus B_3 \oplus B_4 \oplus B_5 $$
such that $A_i \subset V_i$ and $B_i \subset V_i$ for $i = 1,...,5$. 

In other words, we can set $V_1 = D$, $V_2 = E$, $V_3 = F$, etc., and relabel the decompositions
$$A = I \oplus E_1 \oplus G_1 \oplus F \oplus H_1 \oplus H_2$$ 
$$ B = I \oplus E_2  \oplus H_1 \oplus F  \oplus G_1 \oplus G_2$$ 
by setting as $A_i$ as the sum of those summands which lie in $V_i$, i.e. $A_1 = I$, $A_2 = E_1$, $A_3 = F$, $A_4 = G_1$, $A_5 = H_1 \oplus H_2$, and analogously so for $B$.  

Note that for each $i \in \{1,...,5\}$ the subspaces $A_i$ and $B_i$ form a coisotropic pair in $V_i$ of a particularly simple form, each member of the pair being either the entire subspace $V_i$ or a lagrangian subspace therein. Indeed, $A_1 = B_1 = I$ are the same lagrangian subspace of $V_1$, $A_2 = E_1$ and $B_2 = E_2$ form a lagrangian pair in $V_2$, $A_3 = B_3 = F = V_3$, $A_4$ is a lagrangian subspace of $B_4 = G =V_4$, and finally $A_5 = H = V_5$ and $B_5 = H_1$ is lagrangian in this space. We introduce notation for these particularly simple cases of coisotropic pairs. 

\begin{dfn}\label{basic types} A coisotropic pair $(A,B)$ in a symplectic space $V$ is of \textbf{elementary type} if it is one of the following types: 
\begin{description}
\item[$\lambda$:] $A$ and $B$ are lagrangian subspaces, and $A = B$
\item[$\delta$:] $A$ and $B$ are lagrangian subspaces, and $A \cap B = 0$
\item[$\sigma$:] $A = B = V$, i.e. $A$ and $B$ are symplectic
\item[$\mu_B$:] $B = V$ and $A$ is a lagrangian subspace
\item[$\mu_A$:] $A = V$ and $B$ is a lagrangian subspace
\end{description}
We will consider these types ordered as listed and also call them $\tau_1$ through $\tau_5$. 
\end{dfn}

The cases when a coisotropic subspace $C \subset V$ is the entire space or is lagrangian are the two extreme cases of a coisotropic subspace in the sense that they correspond respectively to when $C^{\omega} = 0$ or when $C^{\omega}$ is as large as possible, i.e. $C^{\omega} = C$. The basic types listed above cover all the scenarios when two coisotropics $A$ and $B$ are given by either of these two extremes, except for the possible scenario when $A$ and $B$ are two non-identical lagrangians with non-zero intersection. This case, though, can be split into a ``direct sum'' of the cases $\delta$ and $\lambda$, i.e. it is not ``elementary'' as a type of coisotropic pair. To see this, assume that $A$ and $B$ are such, and let $\tilde A$ and $\tilde B$ be complements of $A \cap B$ in $A$ and $B$ respectively (in particular $\tilde A \cap \tilde B = 0$). Set $W = A + B$ and note that $W^{\omega} = A^{\omega} \cap B^{\omega} = A \cap B \subset W$ because $A$ and $B$ are lagrangian. The subspace $\tilde V := \tilde A \oplus \tilde B$ is such that $\tilde V \oplus (A \cap B) = W$, hence by Lemma \ref{Witt-Artin} it is symplectic and 
$$ V= \tilde V  \oplus \tilde V ^{\omega}$$
with $A \cap B$ as a lagrangian subspace of $\tilde V ^{\omega}$. With respect to this decomposition of $V$, the coisotropics $A$ and $B$ decompose as $A = \tilde A \oplus A \cap B$ and $B = \tilde B \oplus A \cap B$, where $\tilde A$ and $\tilde B$ are a lagrangian pair in $\tilde V$, i.e. a coisotropic pair of type $\delta$, whereas $A \cap B$, seen as the component of both $A$ and $B$ in $\tilde V^{\omega}$, represents a coisotropic pair in $\tilde V^{\omega}$ of the type $\lambda$.

In the following we make more precise the sense in which a coisotropic pair is the direct sum of smaller coisotropic pairs and in which way the elementary types defined above are indeed elementary.  

\begin{dfn}
Given an $\omega$-orthogonal decomposition of $V$ into a finite number $m \in \mathbb{N}$ of symplectic subspaces
$$ V = \bigoplus_{i=1}^m V_i $$
and given subspaces $A_i$,$B_i \subset V_i$ forming a coisotropic pair in $V_i$ for each $i \in \{1,...,m\}$, we say that $(A,B)$ is the direct sum of the coisotropic pairs $(A_i, B_i)$ if
$$ A = \bigoplus_{i=1}^m A_i \quad \quad \text{and} \quad \quad B = \bigoplus_{i=1}^m B_i $$
Such a direct sum decomposition will be denoted 
$$ (A,B) = \bigoplus_{i}^m (A_i,B_i) $$
\end{dfn}

\begin{dfn} A coisotropic pair $(A,B)$ in $V$ is called \textbf{elementary} if there exists no such direct sum decomposition of $(A,B)$ except as a direct sum of coisotropic pairs of only one of the elementary types $\lambda$, $\delta$, $\sigma$, $\mu_B$ or $\mu_A$. 
\end{dfn}

\begin{prop}\label{elementary types are elementary}
The elementary types $\lambda$, $\delta$, $\sigma$, $\mu_B$ and $\mu_A$ are elementary according to the above definition. 
\end{prop}

\pf Assume that $(A,B)$ is a coisotropic pair of some elementary type 
$$\tau \in \{\lambda, \delta, \sigma, \mu_B, \mu_A\}$$
and that 
$$ (A,B) = \bigoplus_{i}^m (A_i,B_i) $$
is a direct sum decomposition into coisotropics, subordinate to an $\omega$-orthogonal decomposition $V  = \bigoplus V_i$ into symplectic subspaces, i.e. such that $A_i \subset V_i$ and $B_i \subset V_i$ for each $i$. We need to show that each coisotropic pair $(A_i,B_i)$ in $V_i$ is of type $\tau$. Because $\tau$ is an elementary type, $A$ is either equal to $V$ or is lagrangian in $V$. If $A = V$, then $A_i = V_i \ \forall i$ for dimension reasons. If $A$ is lagrangian, it is in particular isotropic, and hence each $A_i$ is isotropic in $V$ because $A_i \subset A$. Because $A_i \subset V_i$,  we have $\omega = \omega_{V_i}$ on $A_i$, so $A_i$ is also isotropic in $V_i$. Since $A_i$ is assumed coisotropic in $V_i$, it follows that $A_i$ is lagrangian in $V_i$.  By the same arguments, if $B = V$ then $B_i = V_i \ \forall i$, or if $B$ is lagrangian in $V$ then $B_i$ is lagrangian in $V_i \ \forall i$. It is now clear that if $\tau = \sigma$, then $A_i = B_i = V_i$ for all $i$, so the summand pairs $(A_i,B_i)$ are all also of type $\sigma$. If $\tau =\delta$, then all the $A_i$ and $B_i$ are lagrangian subspaces in their respective $V_i$, and $A \cap B = 0$ implies that $A_i \cap B_i = 0$ for all $i$, so each pair $(A_i,B_i)$ is also of type $\delta$. If $\tau = \lambda$, then similarly the $A_i$ and $B_i$ are lagrangian in $V_i$. To see that here $A_i = B_i \ \forall i$, consider $v \in A = B \subset V$, which has a unique decomposition $v = v_1 + ... + v_m$ with $v_i \in V_i$ for each $i$. Because $v\in A$ and $v \in B$, $v$ also has such unique decompositions with respect to $A = \bigoplus A_i$ and $B = \bigoplus B_i$, but because $A_i,B_i \subset V_i$ for each $i$, these decompositions must coincide with the above decomposition. Hence $v_i \in A_i \cap B_i$ for each $i$. In particular $A = B \subset \bigoplus (A_i \cap B_i)$, which, for dimension reasons, implies $A_i = B_i$ for all $i$.  So each pair $(A_i,B_i)$ is indeed of type $\lambda$ when $(A,B)$ is. Now assume $\tau = \mu_B$. For each $i$, $A_i$ is lagrangian in $V_i$ and $B_i = V_i$, so $(A_i,B_i)$ is also of type $\mu_B$. The case for $\mu_A$ is the same, but with the roles of $A$ and $B$ reversed. \qed

\begin{cor} If a coisotropic pair $(A,B)$ has a direct sum decomposition 
$$ (A,B) = \bigoplus_{i}^m (A_i,B_i) $$
where every coisotropic pair $(A_i,B_i)$ is of the same elementary type, then $(A,B)$ is elementary and of that type. 
\end{cor}

\pf
It suffices to show that $(A,B)$ is of the same type as its summands, since by Proposition \ref{elementary types are elementary} it is then elementary. If the elementary type of the summands is such that the $A_i$ are all lagrangian subspaces of the $V_i$, then the $A_i$ are isotropic subspaces of $V$ and hence their $\omega$-orthogonal sum $A = \bigoplus A_i$ will also be isotropic. Since each $A_i$ has half the dimension of $V_i$, $A$ will have half the dimension of $V$, i.e. it is lagrangian. If on the other hand the elementary type in question is such that  $A_i = V_i$ for each $i$, then clearly $A = V$. The same arguments apply to $B$. Thus the coisotropic pair $(A,B)$ is such that $A$ and $B$ are each either lagrangian or all of $V$ in the same way that their summands $A_i$ and $B_i$ are. It remains only to be sure that when $A$ and $B$ are both lagrangian, they are either identical or such that $A\cap B = 0$, according to whether $A_i = B_i \ \forall i$ or $A_i \cap B_i = 0  \ \forall i$. If $A_i = B_i \ \forall i$ then clearly $A = B$. Assume $A_i \cap B_i = 0 \ \forall i$ and let $v \in A \cap B$. We have a unique decomposition $v = v_1 + ... + v_m$ with $v_i \in V_i$ for all $i$, and because $A = \bigoplus A_i$ and $B = \bigoplus B_i$  are direct sum decompositions subordinate to $V = \bigoplus V_i$, each $v_i$ lies in $A_i \cap B_i = 0$. Hence $v = 0$, and we conclude that $A \cap B = 0$ when $A_i \cap B_i = 0 \ \forall i$.\qed

Proposition \ref{elementary types are elementary} guarantees that the five elementary types of coisotropic pairs are independent of one another in the sense that one cannot express any one of them as a sum of the others. The proof of Proposition \ref{main prop} showed that these basic types are also ``generating'' in the sense that any coisotropic pair decomposes into a direct sum of such elementary types. The corollary implies that one can simplify any direct sum decomposition of a coisotropic pair so that it has only five summands, these summands being of one each of the elementary types. We will call any such five part decomposition an \textbf{\emph{elementary decomposition}}. The following shows that elementary decompositions give a set of invariants for a coisotropic pair $(A,B)$ which are equivalent to the original invariants we associated to such a pair. 

\begin{prop}\label{invariants are equivalent} Let $(A,B)$ be a coisotropic pair in $V$ and let
$$ (A,B) = \bigoplus_{i}^5 (A_i,B_i) $$
be an elementary decomposition subordinate to an $\omega$-orthogonal decomposition $$V= \bigoplus_{i=1}^5 V_i$$ 
ordered such that $(A_i,B_i)$ is of type $\tau_i \in \{ \lambda, \delta, \sigma, \mu_B, \mu_A \}$. Set $n_i := \frac{1}{2} \dim V_i$. Then the 5-tuple 
$$\textbf{n} := (n_1,..., n_5)$$ 
gives a set of invariants (call them \textbf{elementary invariants}) which are equivalent to the canonical invariants 
$$\textbf{k} := (\dim (A^{\omega} \cap B^{\omega}, \ \dim A^{\omega}, \ \dim B^{\omega}, \ \tiny{\frac{1}{2}} \dim V, \  \dim (A^{\omega} \cap B) )$$

\end{prop}

\pf Consider $\textbf{n} = (n_1,...,n_5)$ as a coordinate in the space $\mathcal{N} := \mathbb{Z}_{\scriptscriptstyle{ \geq 0}}^5$ of all possible 5-tuples of elementary invariants (each $V_i$ is symplectic, hence of even dimension), and let $\mathcal{K}$ denote the space of all possible sets of canonical invariants $\textbf{k} = (k_1,...,k_5)$. 

Fix a coisotropic pair $(A,B)$ and fix also an elementary decomposition of this pair, with $A = A_1 \oplus ... \oplus A_5$ and $B = B_1 \oplus ... \oplus B_5$. This gives a 5-tuple $\textbf{n}$. From this $\textbf{n}$ we can obtain the canonical invariants $\textbf{k}$ associated to $(A,B)$ as follows. 

Clearly one has
$$ k_4 =  \tiny{\frac{1}{2}} \dim V = n_1 + n_2 + n_3 + n_4 + n_5 $$
For the remaining invariants, we claim that
$$k_1 = \dim (A^{\omega} \cap B^{\omega}) = n_1 $$
$$ k_2 = \dim A^{\omega} = n_1 + n_2 +  n_4  $$
$$ k_3 = \dim B^{\omega} = n_1 +  n_2 +  n_5 $$
and
$$k_5 = \dim (A^{\omega} \cap B) = n_1 + n_4 $$
To see this, we show 
$$A^{\omega} = A_1 \oplus A_2 \oplus A_4, \quad \quad B^{\omega} =  B_1 \oplus B_2 \oplus B_5$$ 
$$A^{\omega} \cap B = A_1 \oplus A_4, \quad \text{and}  \quad A^{\omega} \cap B^{\omega} = A_1$$
which gives the above formulae for $k_1, k_2,k_3$ and $k_5$ directly. 

For any $a \in A$ we have the decomposition $a = a_1 + ... + a_5$ with $a_i \in A_i$, and for $\tilde a$ also in $A$
\begin{equation}\label{orthogonality splitting}
\omega (a,\tilde a) = \omega_{V_1}(a_1,\tilde a_1) + ... + \omega_{V_5}(a_5,\tilde a_5) = \omega_{V_3}(a_3,\tilde a_3) + \omega_{V_5}(a_5,\tilde a_5)
\end{equation}
because $A_1, A_2$ and $A_4$ are lagrangian in their respective $V_i$. If $\tilde a$ is in $A^{\omega}$, then choosing $a$ as any element  in $A_3$ we find $0 = \omega_{V_3}(a, \tilde a_3)$ and hence $\tilde a_3 \in A_3^{\omega_{V_3}} = 0$, since $A_3$ is symplectic in $V_3$. Similarly one finds $\tilde a_5 = 0$, so $\tilde a \in A_1 \oplus A_2 \oplus A_4$, which shows $A^{\omega} \subset A_1 \oplus A_2 \oplus A_4$. The opposite inclusion $A^{\omega} \supset A_1 \oplus A_2 \oplus A_4$ follows from (\ref{orthogonality splitting}) as well, since for $\tilde a \in A_1 \oplus A_2 \oplus A_4$ and any $a \in A$ we find $\omega (a,\tilde a) = 0$. Arguing analogously one also shows $B^{\omega} = B_1 \oplus B_2 \oplus B_5$.

For the equalities  $A^{\omega} \cap B = A_1 \oplus A_4$ and $A^{\omega} \cap B^{\omega} = A_1$ we use the fact that if $v$ is in $A^{\omega} \cap B$ or $A^{\omega} \cap B^{\omega}$, then in particular $v$ is in $A \cap B$ and hence has a unique decomposition $v = v_1 + ... + v_5$ with $v_i \in A_i \cap B_i \ \forall i$. 

If $v \in A^{\omega} \cap B$, then $v \in A^{\omega} = A_1 \oplus A_2 \oplus A_4$ implies $v_3 = v_5 = 0$. Also, $v_2 \in A_2 \cap B_2 = 0$. Thus $v \in A_1 \oplus A_4$ and $A^{\omega} \cap B \subset A_1 \oplus A_4$ holds. On the other hand, because $A_1 = B_1$ and $A_4 \subset B_4 = V_4$, we have $A_1 \oplus A_4 \subset A^{\omega} \cap B$. 

If $v \in A^{\omega} \cap B^{\omega}$, then not only are $v_3, v_5$ and $v_3$ zero because $A^{\omega} \cap B^{\omega} \subset A^{\omega} \cap B$, but also $v_4 = 0$, because $B^{\omega} = B_1 \oplus B_2 \oplus B_5$ does not contain non-zero summands in $B_4$. Thus $A^{\omega} \cap B \subset A_1$. The opposite in inclusion holds since $A_1 = B_1$ is a summand in the decompositions of both $A^{\omega}$ and $B^{\omega}$.

The equations above describing the $k_i$ in terms of the $n_i$ define a linear map $M: \mathcal{N} \longrightarrow \mathcal{K}$, representable by matrix multiplication with the matrix
$$
M =
\left( \begin{array}{ccccc}
1 & 0 & 0 & 0 & 0 \\
1 & 1 & 0 & 1 & 0 \\
1 & 1 & 0 & 0 & 1 \\
1 & 1 & 1 & 1 & 1 \\
1 & 0 & 0 & 1 & 0 
\end{array} \right)
$$
which is non-singular ($\det M = 1$). Hence $M$ defines an injective map, which means in particular that the numbers $\textbf{n} = (\frac{1}{2} \dim V_1,...,\frac{1}{2} \dim V_5)$ which we associate to an elementary decomposition of a coisotropic pair $(A,B)$ do not depend on the particular elementary decomposition but only depend on the pair $(A,B)$. In other words, $\textbf{n}$ does in fact define a set of invariants for $(A,B)$. The map $M$ is also surjective. Any $\textbf{k} \in \mathcal{K}$ is, by definition, realizable by some coisotropic pair $(A,B)$ and by the proof of Propostion \ref{main prop} this pair has an elementary decomposition; by the above, the invariants $\textbf{n}$ associated to this decomposition are mapped under $M$ to $\textbf{k}$. \qed

To compute the elementary invariants from the canonical invariants one can simply use the inverse of the mapping $M : \textbf{n} \mapsto \textbf{k}$,  
$$
M^{-1} =
\left( \begin{array}{ccccc}
1 & 0 & 0 & 0 & 0 \\
0 & 1 & 0 & 0 & -1 \\
1 & 0 & -1 & 1 & -1 \\
-1 & 0 & 0 & 0 & 1 \\
-1 & -1 & 1 & 0 & 1 
\end{array} \right)
$$
which gives the linear equations for the $n_i$ in terms of the $k_i$:
\begin{align}
n_1 & = k_1\label{n_1 equation} \\
n_2 & = k_2 - k_5 \\
n_3 & = k_1 - k_3 + k_4 - k_5 \\
n_4 & = -k_1 + k_5 \\
n_5 & = -k_1 - k_2 + k_3 + k_5 \label{n_5 equation}
\end{align}
Note that we already nearly explicitly computed these equations in the proof of Proposition \ref{main prop}. 

\begin{cor}\label{ineq only constraint}
The canonical invariants $(k_1,...,k_5)$ are subject only to the five inequalities 
\begin{equation*}
0 \leq k_1 \leq k_5 \leq k_2 \quad \quad k_1 + k_2 \leq k_3 + k_5  \leq k_1 +  k_4 
\end{equation*} 
\end{cor}

\pf
That the $k_i$ must satisfy these inequalities follows from the linear equations (\ref{n_1 equation}) though (\ref{n_5 equation}) for the $n_i$ in terms of the $k_i$ and the fact that $n_i \geq 0 \ \forall i$. The equation for $n_1$ implies $0 \leq k_1$, the equation for $n_2$ gives $k_5 \leq k_2$, the one for $n_3$ gives $k_3 + k_5 \leq k_1 + k_4$, and the inequalities $k_1 \leq k_5$ and $k_1 + k_2 \leq k_3 + k_5$ follow from the equations for $n_4$ and $n_5$. 

To see that these inequalities are the only constraints on the $k_i$, let $\textbf{k} = (k_1,...,k_5)$ be an arbitrary 5-tuple of integers subject only to the above inequalities. We need to show that $\textbf{k}$ is in  $\mathcal{K}$, the set of canonical invariants realizable by a coisotropic pair, which is the image of $M$. In other words we must find a 5-tuple of non-negative integers $\textbf{n} = (n_1,...,n_5)$ such that $M \cdot \textbf{n} = \textbf{k}$, i.e. which solve the linear equations
\begin{align*}
k_1 & = n_1\\
k_2 & = n_1 + n_2 + n_4 \\
k_3 & = n_1 + n_2 + n_5  \\
k_4 & = n_1 + n_2 + n_3 + n_4 + n_5 \\
k_5 & = n_1 + n_4 
\end{align*}
For $k_1 \geq 0$ we choose $n_1 = k_1$ and for $k_5 \geq k_1$ we can always choose $n_4 \geq 0$ such that $k_5 = k_1 + n_4 = n_1 + n_4$. Next, because $k_2 \geq k_5 = n_1 + n_4$, we can choose $n_2 \geq 0$ such that $k_4 = k_5 + n_2 = n_1 + n_2 + n_4$. Thus far $n_1, n_2$ and $n_4$ are fixed and the equations for $k_1, k_2$ and $k_5$ solved. For $k_3$ we have $k_3 \geq k_1 + k_2 - k_5 = n_1 + n_2$, so $n_5$ can be chosen such that $k_3 = n_1 + n_2 + n_5$. Finally, for $k_4 \geq k_3 + k_5 - k_1 = n_1 + n_2 + n_4 + n_5$, an integer $n_3 \geq 0$ is still free to be chosen such that $k_4 = k_3 + k_5 - k_1 + n_3 = n_1 + n_2 + n_3 + n_4 + n_5$ as desired. 
\qed

Using the elementary invariants one can easily construct a normal form $(A_0, B_0)$ for a coisotropic pair $(A,B)$, $\text{i.e.}$ a standardized representative of the equivalence class of $(A,B)$. Let $\textbf{n} = \{n_1,...,n_5\}$ be the elementary invariants of $(A,B)$. We choose $\mathbb{R}^{2n_1} \oplus ... \oplus \mathbb{R}^{2n_5}$ as our model space, equip each summand with the standard symplectic form $\Omega_i$ represented by the $2n_i \times 2n_i$ matrix
$$
\left(  
\begin{array}{cc}
0 & \mathds{1} \\
-\mathds{1} & 0
\end{array}
\right) $$
and give the whole space the direct sum symplectic form $\Omega_1 \oplus ... \oplus \Omega_5$.  Let 
$$(q_1^{i},...,q_{n_i}^{i},p_1^{i},...,p_{n_i}^{i})$$ 
denote the standard coordinates on $\mathbb{R}^{2n_i}$ and denote 
$$Q^{n_i} = \text{span} \{q_1^{i},...,q_{n_i}^{i} \} \text{ and } \ P^{n_i} = \text{span} \{ p_1^{i},...,p_{n_i}^{i} \} $$
Then 
$$ A_0 := Q^{n_1} \oplus Q^{n_2} \oplus \mathbb{R}^{2n_3} \oplus Q^{n_4} \oplus \mathbb{R}^{2n_5}$$
$$ B_0 := Q^{n_1} \oplus P^{n_2} \oplus \mathbb{R}^{2n_3} \oplus \mathbb{R}^{2n_4} \oplus Q^{n_5}$$
defines a normal form for $(A,B)$. By construction $(A_0,B_0)$ is a coisotropic pair such that the elementary invariants of $(A_0,B_0)$ and $(A,B)$ match. Indeed the very definition of $(A_0,B_0)$ gives an elementary decomposition with appropriate dimensions: $(Q^{n_1},Q^{n_1})$ is a coisotropic pair of elementary type $\lambda$ in $\mathbb{R}^{2n_1}$, $(Q^{n_2},P^{n_2})$ a pair of type $\delta$ in $\mathbb{R}^{2n_2}$, and so on. From Proposition \ref{invariants are equivalent} we know that the canonical invariants of $(A,B)$ and $(A_0,B_0)$ match because their elementary invariants do, and by Proposition \ref{main prop} this means that $(A,B) \sim (A_0,B_0)$.

\end{document}